\documentclass[11pt]{article}
\usepackage{amsthm}
\usepackage{amsfonts}
\usepackage[T1]{fontenc}
\usepackage[latin1]{inputenc}
\usepackage{epsfig}

\usepackage[lflt]{floatflt}
\usepackage{graphicx}
\usepackage{amsmath}
\usepackage{amssymb}
\usepackage{fancyhdr}
\usepackage{color}
\usepackage[affil-sl]{authblk}
\usepackage{cancel}
\usepackage{subfigure}
\usepackage[english]{babel}
\usepackage{fancyhdr}
\usepackage{makeidx}
\usepackage{latexsym}

\newtheorem{theorem}{Theorem}

\makeindex \textwidth=17cm \textheight=25cm \oddsidemargin=-0.5cm
\newcommand{\bieq}{\begin{equation}}
\newcommand{\eneq}{\end{equation}}

\newcommand{\bbR}{{\mathbb R}}

\newcommand{\bbC}{{\mathbb C}}

\newcommand{\al}{\alpha}
\newcommand{\la}{\lambda}

\newcommand{\be}{\beta}
\newcommand{\G}{\Gamma}

\begin{document}
\begin{center}
\LARGE
\textbf{A new equivalence of Stefan's problems for the Time-Fractional-Diffusion Equation.}
\end{center}
\medskip
\begin{center}
\normalsize
 Sabrina Roscani$^1$ and Eduardo A. Santillan Marcus$^2$\\
\medskip
\small
$^{1,2}$Departamento de Matem\'{a}tica,  FCEIA, Universidad Nacional de Rosario, Pellegrini 250, Rosario, Argentina \\
\textcolor{blue}{ sabrina@fceia.unr.edu.ar, edus@fceia.unr.edu.ar} \\

$^2$ Departamento de Matem\'atica,
FCE, Universidad Austral, Paraguay 1950, Rosario, Argentina,\\
\textcolor{blue}{esantillan@austral.edu.ar}\\

$^1$ CONICET, Argentina.

\end{center}

 \bigskip \medskip
\textbf{Note:} This paper was already accepted to be published in the in the journal "Fractional Calculus and Applied Analysis".

%%%% Abstract %%%%%%%%%%%%%%%%%%%%%%%%%
 \medskip
\small
\noindent \textbf{Abstract: }

A fractionary Stefan problem with a boundary convective condition is solved, where  the fractional derivative of order $ \al \in (0,1) $
is taken in the Caputo sense. Then an equivalence with other two fractional Stefan problems (the first one with a constant condition
on $ x = 0 $ and the second with a flux condition)is proved and  the convergence to the classical solutions is
analyzed when  $ \al \nearrow $ 1 recovering the heat equation with
its respective Stefan's condition.
\medskip
\normalsize


 \smallskip


 \section{Introduction}\label{sec:1}

\setcounter{section}{1}
\setcounter{equation}{0}\setcounter{theorem}{0}

\noindent In 1695 L'Hopital inquired of Leibnitz, the father of the
concept of the classical differentiation, what meaning could be
ascribed to the derivative of order $\frac{1}{2}$. Leibnitz replied
prophetically: ``[...] this is an apparent paradox from which, one
day, useful consequences will be drawn.''

\noindent From 1819, mathematicians as Lacroix, Abel, Liouville,
Riemann and later Gr\"{u}nwald and Letnikov attempted to establish a
definition of fractional derivative.

\noindent  We use here the definition introduced by Caputo in 1967,
and we will call it \textsl {fractional derivative in the Caputo's
sense}, which is given by

$$\,_{a} D^{\alpha}f(t)=D^{\alpha}f(t)=\frac{1}{\Gamma(n-\alpha)}\int^{t}_{a}(t-\tau)^{n-\alpha-1} f^{(n)}(\tau)d\tau $$
where $\al >0$ is the order of derivation, $n=\left\lceil \al
\right\rceil$ and  $f $ is a differentiable function up to order $n$
in $\left[a,b\right]$. To simplify notation, we use from here the
notation $D^{\alpha}$ for the fractional derivative in the Caputo's
sense.

The one-dimensional heat equation has become the paradigm for the
all-embracing study of parabolic partial differential equations,
linear and nonlinear. Cannon \cite{Cannon} did a methodical
development of a variety of aspects of this paradigm. Of particular
interest are the discussions on the one-phase Stefan problem, one of
the simplest examples of a free-boundary-value problem for the heat
equation (see  Datzeff \cite{dat}). In mathematics and its
applications, particularly related to phase transitions in matter, a
Stefan problem is a particular kind of boundary value problem for a
partial differential equation, adapted to the case in which a phase
boundary can move with the time. The classical Stefan problem aims
to describe the temperature distribution in a homogeneous medium
undergoing a phase change, for example ice passing to water: this is
accomplished by solving the heat equation imposing the initial
temperature distribution on the whole medium, and a particular
boundary condition, the Stefan condition, on the evolving boundary
between its two phases. Note that in the one-dimensional case this
evolving boundary is an unknown curve: hence, the Stefan problems
are examples of free boundary problems. A large bibliography on free
and moving boundary problems for the heat-diffusion equation was
given in Tarzia \cite{Tar2}.

\smallskip

In this paper, we deal with three one-phase Stefan's problems with
time fractional diffusion equation, obtained from the standard
diffusion equation by replacing the first order time-derivative by a
fractional derivative of order $\alpha > 0 $ in the Caputo sense:
$$  \, D^{\alpha}u(x,t)=\lambda^2\dfrac{\partial^2 u }
{\partial x^2}(x,t), \quad  -\infty<x<\infty, \ t>0, \ 0<\al<1, $$
 and the Stefan condition $\dfrac{d s(t)}{dt}=k u_x(s(t),t), \, t>0,$ by the fractional Stefan condition
$$  \, D^{\alpha}s(t)= k u_x(s(t),t),\quad t>0. $$

 The fractional diffusion equation has been treated by a number of authors (see  \cite{FM-AnalPropAndAplOfTheWFunc}, \cite{FM-TheFundamentalSolution},
\cite{Kilbas},\cite{Luchko},
 \cite{Podlubny}) and, among the several applications that
have been studied, Mainardi \cite{FM-libro} studied the application
to the theory of linear viscoelasticity.\\

Fractional moving boundary problems have gained in recent years
great interest for the applications and these contributions have a
potential significant impact because exact solutions and equivalence
of different problems are provided. To have a complete review of the
results in this field, as a reference for further studies to know
the mathematical results obtained in literature and their
applications, see \cite{Atkinson}, \cite{Kushwaha},
\cite{LiWangZhao}, \cite{Liu-Xu}, \cite{Vogl}, \cite{Voller}. An
interesting physical meaning of the fractional Stefan's problems is
discussed in \cite{Garra}.\\

\section{Some previous results}\label{sec:2}

\setcounter{section}{2}
\setcounter{equation}{0}\setcounter{theorem}{0} Let us consider the
following problems

\begin{equation}{\label{St1}}
\left\{\begin{array}{lll}
          D^{\al} u(x,t)=\lambda^2\dfrac{\partial^2u }{\partial x^2}(x,t), &   0<x<s(t), \,  t>0, \,  0<\al<1 , \, \, \la>0,\\

          u(0,t)=B,   &  t>0, \quad  B>0 \, \text{ constant}, \\

          u(s(t),t)=C<B, & t>0,\\
          D^{\al}s(t)=-k u_x(s(t),t), & t>0, \quad  k>0, \, \text{ constant}
          \\
          s(0)=0,                                   \end{array}\right.\end{equation}

\noindent and
\begin{equation}{\label{St2}}
\left\{\begin{array}{lll}
          D^{\al} u(x,t)=\lambda^2\dfrac{\partial^2u }{\partial x^2}(x,t), &  & 0<x<s(t), \, t>0, \, 0<\al<1, \, \, \la>0, \\

          u_x(0,t)=-\frac{q}{t^{\al/2}}, \quad  \quad &   & t>0, \quad q>0 \, \text{ constant},\\

          u(s(t),t)=C, \quad  \quad & & t>0,\\
          D^{\al}s(t)=-k u_x(s(t),t),\quad  \quad & & t>0,
          \\
          s(0)=0.                                   \end{array}\right.\end{equation}

\medskip

\noindent A pair $\{u,s\} $ is a solution of the problem
$(\ref{St1})$ (or $(\ref{St2})$) if

\begin{enumerate}
    \item $u$ and $s$ satisfy (\ref{St1}) (or (\ref{St2})),
    \item $u_{xx} $ and  $D^{\al} u$ are continuous for $0<x< s(t)$, $0<t<T$,
    \item $u$ and $u_x$ are continuous for $0\leq x\leq s(t)$,
    $0<t<T$,
    \item $0\leq \underset{x,t\rightarrow 0^+}{\liminf}\,u(x,t)\leq \underset{x,t\rightarrow 0^+}{\limsup }\,
    u(x,t)<+\infty$,
    \item $s$ is continuously differentiable in $[0,T)$ and $\dfrac{\dot{s}(\tau)}{(t-\tau)^\al}$ $\in L^1(0,t)$ $\forall t\in (0,T)$.
    \end{enumerate}
There are two important functions involved in the solution of this kind of problems (see \cite{FM-theM-Wright}):\\

\noindent The Wright function
$$W(z,\al,\be)=\sum^{\infty}_{n=0}\frac{z^n}{n! \G\left( \al n+\be \right)},\, z\in \bbC\, \, , \, \al>-1,$$
and  a particular case, the Mainardi function
 \begin{equation} \label{M_nu} M_\nu (z)= W(-z,-\nu,1-\nu)=\sum^{\infty}_{n=0}\frac{(-z)^n}{n! \G\left( -\nu n+ 1-\nu \right)}.  \end{equation}

These two problems were solved in \cite{RoSa} and its solutions are given by\\
\begin{equation}{\label{SOL-St1}}
\left\{\begin{array}{lll}
           u_1(x,t)=B+\frac{C-B}{1-W\left(-\tilde{\xi},-\frac{\al}{2},1\right)}
           \left[1-W\left(-\frac{x}{\la t^{\al/2}},-\frac{\al}{2},1\right)\right],\\
s_1(t)=\la \tilde{\xi} t^{\al/2}, \\
          \text{where } \tilde{\xi} \text{ is the unique solution to the equation}  \\

          H(\xi)=-\frac{k}{\la^2}
 \frac{\G(1-\frac{\al}{2})}{\G(1+\frac{\al}{2})}(C-B),                                   \end{array}\right.\end{equation}

where
\begin{equation}\label{St1-H(xi)}
H(\xi)=\xi\left[1-W\left(-\xi,-\frac{\al}{2},1\right)\right]
 \frac{1}{M_{\al/2}\left(\xi\right)},
 \end{equation}

and
\begin{equation}{\label{SOL-St2}}
\left\{\begin{array}{lll}
           \begin{array}{ll}
             u_2(x,t)= & C+q\la\G\left(1-\frac{\al}{2}\right)\left[ 1-W\left(-\tilde{\mu},-\frac{\al}{2},1\right)\right]- \\
              & -q\la\G\left(1-\frac{\al}{2}\right)
           \left[ 1-W\left(-\frac{x}{\la
           t^{\al/2}},-\frac{\al}{2},1\right)\right],
           \end{array}
           \\
s_2(t)=\la \tilde{\mu} t^{\al/2}, \\
          \text{where } \tilde{\mu} \text{ is the unique solution to  equation}  \\

          J(\mu)=\frac{kq}{\la}
 \frac{\G\left(1-\frac{\al}{2}\right)^2}{\G(\frac{\al}{2}+1)},        \end{array}\right.\end{equation}

 where

 \begin{equation}\label{St2-J(xi)}
J(\mu)=\mu \frac{1}{ M_{\al/2}(\mu)}.
 \end{equation}

As we can see, the Mainardi function (\ref{M_nu}), and the "fractional \emph{erf} function" ($ 1-W\left(-x,-\frac{\al}{2},1\right) $) are essential in this study. \\

As it seen in \cite{RoSa}, if  $\,  0<\al<1$ and $x \in \bbR^+$, the
Mainardi function $M_{\al/2}(x) $ is a decreasing positive function
, and  $1-W\left(-x,-\frac{\al}{2},1\right)$ is a increasing
positive function.

\begin{theorem}Let us consider problems (\ref{St1}) and (\ref{St2}) where
 \begin{enumerate}
 \item the constant $C$ is the same in both problems,
 \item in problem (\ref{St1}) $B=C-q\la\G\left(1-\frac{\al}{2}\right)\left[ 1-W\left(-\tilde{\mu},-\frac{\al}{2},1\right)\right],$ where $\tilde{\mu}$ is the unique solution to
     $J(\mu)=\dfrac{kq}{\la}
 \frac{\G\left(1-\frac{\al}{2}\right)^2}{\G(\frac{\al}{2}+1)}$ and $J$ is defined by
 $(\ref{St2-J(xi)})$.
 \end{enumerate}
 \noindent Then these problems are equivalent.
\end{theorem}

\proof See \cite{RoSa}, page 808.

\section{A Fractional Stefan Problem with a boundary convective condition}
Let us consider now the following problem

\begin{equation}{\label{St3}}
\left\{\begin{array}{lll}
          D^{\al} u(x,t)=\lambda^2\dfrac{\partial^2u }{\partial x^2}(x,t), &   0<x<s(t), \,  t>0, \,  0<\al<1 , \, \, \la>0,\\
          m u_x(0,t)=\frac{h}{t^{\al/2}}(u(0,t)-D),    &  t>0, \quad  D>0 \, \text{ constant}, \\
          u(s(t),t)=C<D, & t>0,\\
          D^{\al}s(t)=-k u_x(s(t),t), & t>0, \quad  k>0, \\
          s(0)=0.                                   \end{array}\right.\end{equation}

\medskip

\noindent Taking into account the results mentioned in \cite{RoSa},
we propose the following solution
\begin{equation}\label{u}
u(x,t)= a + b W\left(-\frac{x}{\la t^{\al/2}},
-\frac{\al}{2},1\right),
\end{equation}
where the constants $a$ and $b$ will be determined.
\\
$$ u(s(t),t)=C \Leftrightarrow a+bW\left(-\frac{s(t)}{\la t^{\al/2}}, -\frac{\al}{2},1\right)=C \quad \forall t.$$
Due to  the strict monotonicity of $u$,this expression is valid for
every $t>0$  only if $s(t)$ is proportional to $\la t^{\al/2}$ ,

\begin{equation}\label{s(t)}
s(t)=\eta \la t^{\al/2}\end{equation}

\noindent then,
\begin{equation} \label{1} C=a+bW\left(-\eta, -\frac{\al}{2},1\right). \end{equation}

$$ u_x(x,t)=-\frac{b}{\la t^{\al/2}}M_{\al/2}\left(\frac{x}{\la t^{\al/2}}\right)\Rightarrow u_x(0,t)= -\frac{b}{\la t^{\al/2}}\frac{1}{\G(1-\al/2)},$$

\noindent and using the convective condition, we have

\begin{equation}
\label{2} -\frac{mb}{\la
t^{\al/2}}\frac{1}{\G(1-\al/2)}=\frac{h}{t^{\al/2}}(a+b-D).
\end{equation}

\noindent From (\ref{1}) and (\ref{2})
$$\begin{cases}
a=D-\left(1+\frac{m}{h\la \G(1-\al/2)}\right)\frac{D-C}{1-W\left(-\eta, -\frac{\al}{2},1\right)+\frac{m}{h\la \G(1-\al/2)}},  \\
b=\frac{D-C}{1-W\left(-\eta, -\frac{\al}{2},1\right)+\frac{m}{h\la
\G(1-\al/2)}}
\end{cases}$$

 \noindent and finally,

\begin{equation} \label{u definitiva}
u(x,t)=D-\frac{(D-C)\left[1-W\left(-\frac{x}{\la
t^{\al/2}},-\frac{\al}{2}, 1\right)+\frac{m}{h\la
\G(1-\al/2)}\right]} {1-W\left(-\eta,-\frac{\al}{2},
1\right)+\frac{m}{h\la \G(1-\al/2)}}.
\end{equation}

\noindent Let us work with the fractional Stefan condition. Taking
into account that

$$ D^{\al}(t^{\be})=\frac{\G(\be+1)}{\G(1+\be-\al)}t^{\be-\al} \quad \text{if } \be>-1 $$

\noindent it follows that

\begin{equation}\label{Deriv al St}
D^{\al}s(t)=D^{\al}(\la \eta t^{\al/2})=\la
\eta\frac{\G(\frac{\al}{2}+1)}{\G(1-\frac{\al}{2})}t^{-\al/2}.
\end{equation}

\noindent On the other hand,

 \begin{equation}\label{u_x(s(t))}
 u_{x}(s(t),t)=-b \frac{1}{\la t^{\al/2}}M_{\al/2}\left(\eta\right)=
 -\frac{D-C}{1-W\left(-\eta, -\frac{\al}{2},1\right)+\frac{m}{h\la \G(1-\al/2)}}\frac{1}{\la t^{\al/2}}M_{\al/2}\left(\eta\right).\end{equation}

\noindent Replacing $(\ref{Deriv al St})$ and $(\ref{u_x(s(t))})$ in
the fractional Stefan condition,

\begin{equation}\label{St-eq para xi}
 \eta\left[1-W\left(-\eta,-\frac{\al}{2},1\right)+\frac{m}{h\la \G(1-\al/2)}\right]
 \frac{1}{M_{\al/2}\left(\eta\right)}=\frac{k}{\la^2}
 \frac{\G(1-\frac{\al}{2})}{\G(1+\frac{\al}{2})}(D-C).
 \end{equation}

\noindent Let us defined the function
\begin{equation}\label{St-K(eta)}
K(\eta)=\eta\left[1-W\left(-\eta,-\frac{\al}{2},1\right)+\frac{m}{h\la
\G(1-\al/2)}\right]
 \frac{1}{M_{\al/2}\left(\eta\right)}.
 \end{equation}
$K$ has the following properties:
\begin{enumerate}
\item $K(0^+)=0$,
\item $K(+\infty)=+\infty$,
\item K is continuous and monotonically increasing.
\end{enumerate}

Because of the asymptotic behavior of the Wright function (see \cite{FM-AnalPropAndAplOfTheWFunc}), it is easy to check properties 1 and 2.\\

\noindent For property 3, we observe that,
$1-W\left(-\eta,-\frac{\al}{2},1\right)$ is a positive and increasing function in $\bbR^+$, and $\frac{1}{M_{\al/2}(\eta)}$ is a positive increasing function.\\

Finally, noting that  $\frac{k}{\la^2}
 \frac{\G(1-\frac{\al}{2})}{\G(1+\frac{\al}{2})}(D-C)>0$, we can affirm that there exists a unique  $\tilde{\eta}$ such that
\begin{equation}\label{K(eta)=algo} K(\tilde{\eta})=\frac{k}{\la^2}
 \frac{\G(1-\frac{\al}{2})}{\G(1+\frac{\al}{2})}(D-C). \end{equation}

\noindent  Then, the solution of problem (\ref{St3}) is given by
\begin{equation}{\label{SOL-PBSt1}}
\left\{\begin{array}{lll}
u_3(x,t)=D-\frac{(D-C)\left[1-W\left(-\frac{x}{\la
t^{\al/2}},-\frac{\al}{2},
1\right)+\frac{m}{h\la\G(1-\al/2)}\right]}
{1-W\left(-\tilde{\eta},-\frac{\al}{2}, 1\right)+\frac{m}{h\la \G(1-\al/2)}},\\
s_3(t)=\la \tilde{\eta} t^{\al/2}, \\
          \text{where } \tilde{\eta} \text{ is the unique solution of}  \\

          K(\eta)=\frac{k}{\la^2}
 \frac{\G(1-\frac{\al}{2})}{\G(1+\frac{\al}{2})}(D-C).                                   \end{array}\right.\end{equation}

 \begin{theorem}Let us consider problems (\ref{St1}) and (\ref{St3}) where
 \begin{enumerate}
 \item the constant $C$ is the same in both problems,
 \item in problem (\ref{St1})
 $B=D-(D-C)\frac{m}{h\la\G(1-\al/2)}\frac{1}
{1-W\left(-\tilde{\eta},-\frac{\al}{2}, 1\right)+\frac{m}{h\la
\G(1-\al/2)}}$, where $\tilde{\eta}$ is the unique solution to
$K(\eta)=\frac{k}{\la^2}
 \frac{\G(1-\frac{\al}{2})}{\G(1+\frac{\al}{2})}(D-C)$ and $K$ is defined by $(\ref{St-K(eta)})$.
 \end{enumerate}
 \noindent Then these problems are equivalent.
 \end{theorem}

\proof \noindent Let us define the following function
$$ B(\xi)=D-(D-C)\frac{m}{h\la\G(1-\al/2)}\frac{1}
{1-W\left(-\xi,-\frac{\al}{2}, 1\right)+\frac{m}{h\la
\G(1-\al/2)}}.$$

\noindent Observe that $B(\tilde{\eta})=B>C$.\\
Now $$H(\xi)=-\frac{k}{\la^2}
 \frac{\G(1-\frac{\al}{2})}{\G(1+\frac{\al}{2})}(C-B(\xi)) \Longleftrightarrow $$

$$ \xi\left[1-W\left(-\xi,-\frac{\al}{2},1\right)\right]
 \frac{1}{M_{\al/2}\left(\xi\right)}=$$
 $$=-\frac{k}{\la^2}\frac{\G(1-\frac{\al}{2})}{1+\frac{\al}{2}}
 \left[ C-D+(D-C)\frac{m}{h\la\G(1-\al/2)}\frac{1}
{1-W\left(-\xi,-\frac{\al}{2}, 1\right)+\frac{m}{h\la \G(1-\al/2)}}
\right]\Longleftrightarrow$$

$$
 \xi\left[1-W\left(-\xi,-\frac{\al}{2},1\right)\right]
 \frac{1}{M_{\al/2}\left(\xi\right)}=$$

 $$=-(D-C)\frac{k}{\la^2}\frac{\G(1-\frac{\al}{2})}{1+\frac{\al}{2}}
  \frac{-\frac{h\la\G(1-\al/2)}{m}\left[1-W\left(-\xi,-\frac{\al}{2}, 1\right)\right]}{\frac{h\la \G(1-\al/2)}{m}\left[1-W\left(-\xi,-\frac{\al}{2}, 1\right)\right]+1} \Longleftrightarrow$$

\begin{equation}\label{K teo eq}
\Longleftrightarrow
K(\xi)=(D-C)\frac{k}{\la^2}\frac{\G(1-\frac{\al}{2})}{1+\frac{\al}{2}}.
 \end{equation}

\noindent Then if $\tilde{\eta}$ is the unique solution of (\ref{K
teo eq}) we have
 $$H(\tilde{\eta})=\frac{k}{\la^2}
 \frac{\G(1-\frac{\al}{2})}{\G(1+\frac{\al}{2})}(C-B).$$
%and that there exists a unique solution of this equation, we can assure that

Due to the uniqueness of solution of $(\ref{K(eta)=algo})$, we can
assure that $ \tilde{\eta}=\tilde{\xi} $ and therefore $ s_1=s_3 $.\\
 Finally, it is easy to verify that $u_1=u_3$.
\endproof

Analogously, we have the following result.

 \begin{theorem}Let us consider problems (\ref{St2}) and (\ref{St3}) where
 \begin{enumerate}
 \item the constant $C$ is the same in both problems,
 \item in problem (\ref{St2})
 $q=\frac{(D-C)k}{m+h\la \G(1-\frac{\al}{2})\left[1-W(-\tilde{\eta},-\frac{\al}{2},1)\right] }$, where $\tilde{\eta}$ is the unique solution to $K(\eta)=\frac{k}{\la^2}
 \frac{\G(1-\frac{\al}{2})}{\G(1+\frac{\al}{2})}(D-C)$ and $K$ is defined by $(\ref{St-K(eta)})$.
 \end{enumerate}
 \noindent Then these problems are equivalent.
 \end{theorem}

\section{Convergence}
We proved in \cite{RoSa} that if $x\in \bbR^+_0$ and $\al \in (0,1)$, then $ \lim\limits_{\al\nearrow 1}\left[1-W\left(-x,-\frac{\al}{2},1\right)\right]= erf\left(\frac{x}{2}\right).$\\

Applying this result to calculate the limit when $\al \nearrow 1$ to
the given solution (\ref{SOL-PBSt1}),  we recover the solution given
by Tarzia in \cite{Tarzia3}:

$$\lim_{\al\nearrow 1}u_3(x,t)=\lim_{\al\nearrow 1}D-\frac{(D-C)\left[1-W\left(-\frac{x}{\la t^{\al/2}},-\frac{\al}{2}, 1\right)+\frac{m}{h\la\G(1-\al/2)}\right]}
{1-W\left(-\tilde{\eta},-\frac{\al}{2}, 1\right)+\frac{m}{h\la
\G(1-\al/2)}}=
$$
$$ =D-\frac{(D-C)\left[erf\left( \frac{x}{2\la \sqrt{t}} \right){x}+\frac{m}{h\la\sqrt{\pi}}\right]}
{erf(\eta/2)+\frac{m}{h\la \sqrt{\pi}}}, $$
%s_3(t)=\la \tilde{\eta} t^{\al/2} = $$
% $$ =B+\frac{C-B}{erf\left(\tilde{\xi}/2\right)}\, erf\left(-\frac{x}{2\la t^{1/2}}\right). $$

$$\lim_{\al\nearrow 1}s_3(t)=\lim_{\al\nearrow 1}\la \tilde{\eta} t^{\al/2}=\la \tilde{\eta}\sqrt{t},
$$
where $ \tilde{\eta}$  is the unique solution to the equation

$$\eta\left[erf(\eta/2)+\frac{m}{h\la \sqrt{\pi}}\right]
 exp(\eta^2/4)=\frac{k}{\la^2}
 \frac{(D-C)}{\sqrt{\pi}},$$

\section{Conclusions}\label{sec:4}
Continuing with our work \cite{RoSa}, we solved a new fractional
Stefan's Problem with a convective boundary condition and then we
proved the equivalence between this problem and  the other two
fractional Stefan's problems presented in the mentioned work.
Finally, we analyzed the convergence when $\al\nearrow 1$, and we
recovered the solution to the classical Stefan's problem with
convective boundary condition.

\section*{Acknowledgements}

\noindent This paper has been sponsored by the Project PICTO AUSTRAL
2008 Nº 73 from {\it Agencia Nacional de Promoci\'{o}n Cient\'{\i}fica y
Tecnol\'{o}gica de la Rep\'{u}blica Argentina} (ANPCyT) and the project
ING349 {\it "Problemas de frontera libre con ecuaciones
diferenciales fraccionarias"}, from Universidad Nacional de Rosario,
Argentina. We appreciate the valuable suggestions by the anonymous
referees which improve the paper.

%%%%%%%%%%%%%%%%%%%%%%%%%%%%%%%%%%%%%%%%%%%%%

%%%%%%%%%% References %%%%%%%%%%%%%%%%%%%%%%%%%%%%%%%%
%%%% arranged in ALPHABETIC ORDER of Authors' Families

%%%% example for a book %%%%%%%%%%%%%
%\bibitem{GasRah}
% G. Gasper, M. Rahman,
% {\it Basic Hypergeometric Series}.
% Cambridge University Press, Cambridge (1990).

%%%% example for article in FCAA journal %%%%%%%%%%%%%%%%%
%\bibitem{Kir}
% V. Kiryakova,
% A brief story about the operators of generalized
%fractional calculus.
% \emph{Fract. Calc. Appl. Anal.} \textbf{11}, No 2 (2008), 201-218.

%%%% example for journal's article %%%%%%%%%%%%%%%%%
%\bibitem{Moak}
% D.S. Moak,
% The $q$-analogue of the Laguerre polynomials.
%{\it J. Math. Anal. Appl.} {\bf 81}, No 1 (1981), 20-47.

%%%% example for a paper in Proceedings %%%%%%%%%%%%%
%\bibitem{Rosbl}
% M. Rosenblum,
% Generalized Hermite polynomials and the Bose-like oscillator
% calculus.
% In: {\it Operator Theory: Advances and Applications},
% Birkh\"auser, Basel (1994), 369-396.
%%%%%%%%%%%%%%%%%%%%%%%%%%%%%%%%%%%%%%%%%%%%%%%%%%%%

%\end{thebibliography} %%%%%%%%%%%%%%%%%%%%%%%%%%%%%%%%

%%%%%%%%%% put authors' addresses here, in \it %%%%%%%%

 \bigskip \smallskip

 \it

 \noindent
%(First) Author's full postal address
$^1$ Departamento de Matem\'{a}tica - ECEN \\
Facultad de Cs. Exactas, Ingenier\'ia y Agrimensura \\
Universidad Nacional de Rosario \\
Av. Pellegrini 250\
(2000) Rosario, ARGENTINA  \\[4pt]
e-mail: sabrinaroscani@gmail.com
%\hfill Received: December 28, 2012
\\[12pt]
 \noindent
% Second Author's address
$^2$ Departamento de Matem\'{a}tica \\
Facultad de Cs. Empresariales\\
Universidad Austral Rosario \\
Paraguay 1950\
(2000) Rosario, ARGENTINA  \\[4pt]
e-mail: edus@fceia.unr.edu.ar

\smallskip

{\bf Correspondence should be sent to the second author.}

\end{document}